\documentclass{amsart}
\usepackage{amsmath,amsthm,amssymb,IMjournal}

\begin{document}
\newtheorem*{theorem}{Theorem}
\newtheorem{corollary}{Corollary}
\newtheorem{lemma}{Lemma}

\def\gp#1{\langle #1 \rangle}

\numberwithin{equation}{section}

\title{Group algebras whose   groups of normalized units have   exponent $4$}

\author{\| Victor  |Bovdi|, Al Ain,
        \|Mohammed |Salim|, Al Ain}



\abstract
We give a full description of  locally finite $2$-groups $G$ such that the normalized group
of units  of the  group algebra $FG$ over a field $F$ of characteristic $2$ has
exponent $4$.
\endabstract



\thanks
Supported by UAEU grants:  UPAR  G00001922 and StartUP G00001889
\endthanks

\section{Introduction and result}\label{sec1}
It is well known that there does  not exist  an effective description of finite
groups of prime square    exponent $p^2$ (not even in the case when the exponent
is $4$). However  Z.~Janko  (see for example \cite{Janko_1, Janko_2, Janko_3}) was able
to characterize these  groups  under certain  additional restrictions  on their
structure. In this way he obtained  interesting  classes  of finite $p$-groups.

Note also that   there is no  effective description of  finite $2$-groups
with    pairwise commuting involutions. On the other hand, the structure of
a locally finite $2$-group $G$ is known when  its normalized group of units
$V(FG)$ of the group algebra $FG$  has the  property that its   involutory
units pairwise commute (see \cite{Bovdi_Dokuchaev}).

There is a similar situation in the case of  powerful $p$-groups.
Despite of extensive current research of this field,
their structure has been  incompletely described.
However, it is possible to determine  \cite{Bovdi_powerful} those cases
when the normalized groups of units  of the group algebras
are powerful $p$-groups.

So  it is a natural question  whether it is  possible to give a description
of those modular group algebras whose groups of  normalized units have
exponent $p^2$. In this note we deal with the case of $p = 2$.

Note that in the case when p is even,    $2$-groups of exponent $4$ appear
in several problems in group theory and in the theory of group ring units
(see for example  \cite{Bovdi_Survey, Caranti, Gupta_Newman, Quick, Vaughan-Lee_Wiegold}).

Our main result is the following.
\begin{theorem}
Let $V(FG)$ be the normalized group of units  of a  group algebra  $FG$ of a
locally finite nonabelian $2$-group $G$ over a field $F$ with  $char(F)=2$.
The group    $V(FG)$  has exponent $4$ if and only if $G=H\times W$,
where  $H$ is a nilpotent group of class $2$,  the Frattini subgroup  of $H$
is central elementary abelian, $|H'|\leq 4$,  and $W$ is an  abelian group
of exponent at most $4$.
\end{theorem}

\section{Preliminaries and  the proof of the Theorem}\label{sec2}

An involution in a group $G$ is an element of order $2$. For any $a,b\in  G$, we denote
$(a, b) = a^{-1}b^{-1}ab$ and $a^b = b^{-1}ab$. Let $D_8$ and $Q_8$ be the dihedral and
quaternion groups of order $8$, respectively.
Define the following groups:
\[
\begin{split}
G_{16}^3=\gp{g,h\mid \quad  g^4=h^2&=1, \quad (g^2,h)=1,\\
\quad  (gh)^3&=hg^3}\cong (C_4\times C_2)\rtimes C_2;\\
G_{16}^4=\gp{g,h\mid \quad  g^4=h^4&=1,  \quad  g^h=g^3}\cong C_4\rtimes C_4;\\
G_{32}^2=\gp{g,h\mid \quad  g^4=h^4&=(gh)^2=1,\\
 (g^2,h)&=(g,h^2)=1}\cong (C_4\times C_2)\rtimes C_4;\\
G_{32}^6=\gp{g,h\mid \quad  g^4=h^4&=(g^3h)^2=1,\\
 (g^2,h)&=(g,h^2)=1}\cong \big((C_4\times C_2)\rtimes C_2\big)\rtimes C_2.\\
\end{split}
\]
For the designation of these groups $G_m^s$ we use their numbers $s$  in the Small Groups
Library of order $m$ in the computer algebra program  GAP \cite{GAP}.

We use freely  the following well known equations (see \cite{Hall_book}, p. 171)
\begin{equation}\label{E:1}
(a,bc) = (a,b)(a,c)(a, b, c), \quad (ab, c) = (a,c )(a, c, b)(b, c).
\end{equation}
\begin{lemma}\label{L:1}
Let $char(F)=2$ and let $H$ be a nonabelian two-generated subgroup of a group $G$. If
$V(FG)$ has exponent $4$, then
$H\in\{D_8,\;  Q_8, \; G_{16}^3,  \; G_{16}^4, \;  G_{32}^2\}$.
Moreover $H'\subseteq \Phi(H)\subseteq \zeta(H)$.
\end{lemma}
\proof Clearly $G$ has exponent $4$ and  any two involutions  in  $G$  either commute or
generate a dihedral group $D_8$ of order $8$. Put $H=\gp{g,h\mid g,h\in G, (g,h)\not=1}$.
Consider the following cases:

{\underline{Case A}}. Let $H=\gp{g,h\mid |g|=4, |h|= 2}$, such that
$H\not\cong D_8$. Then  $x=1+g+h\in V(FG)$ has order $4$ and
\[
\begin{split}
x^4-1=(gh)^2&+(hg)^2+g^3h+ghg^2\\
&+g^2hg +\underline{hg^3}+g^2+hg^2h=0.
\end{split}
\]
Comparing  $hg^3$ with other elements, from the last equation we get  $hg^3=g^2hg$, so
\begin{equation}\label{E:2}
 (h,g^2)=1 \quad\text{and}\quad(gh)^2=(hg)^2.
\end{equation}
It is easy to see $G/\gp{g^2}$ is generated by two involutions $g\gp{g^2}$ and
$h\gp{g^2}$. In the case when  $G/\gp{g^2}$ is abelian we have $(g,h)=g^2$, it follows that
$hgh=g^3$ and $\gp{g,h}\cong {D}_8$, a contradiction. Therefore, $G/\gp{g^2}\cong D_8$. Since $G$ is not  $D_8$,  by (\ref{E:2}) we obtain that   $(gh)^2=(hg)^2\not= 1$.  Consequently,
$(gh)^3=hg^3=(gh)^{-1}$ and
\[
H=\gp{g,h\mid \quad  g^4=h^2=1, (g^2,h)=1, (gh)^4=1 \quad}\cong G_{16}^3.
\]

\underline{Case B}. Let $H=\gp{g,h\mid |g|=|h|=4, |gh|=2}$. Clearly  the unit  $y=1+g+gh\in V(FG)$
has order $4$. Since $ghg=h^3$,  we have that   $y^2=g^2+g^2h+h^3$ and
\[
\begin{split}
y^4-1=gh^3gh&+h^2+\underline{h}+gh^3g\\
&+g^2h^3+h^3g^2+g^2+h^3g^2h=0.
\end{split}
\]
Comparing  $h$ with other elements, from the last equation we obtain   that  only    $h=gh^3g$, so    $(h,g^2)=1$.  Consequently
\[
H=\gp{\quad g,h\mid g^4=h^4=(gh)^2=1,(g^2,h)=1\quad }\cong G_{16}^3.
\]
\underline{Case C}. Let $H=\gp{g,h\mid |g|=|h|=|gh|=4}$ such that  $H\not\cong Q_8$. Then $x=1+g+h\in V(FG)$ has order $4$ and
\begin{equation}\label{E:3}
\begin{split}
x^4-1=(gh)^2+(hg)^2&+\underline{g^3h}+ghg^2+g^2h^2\\
+h^2g^2&+g^2hg+hg^3+h^2gh+gh^3\\
&+h^3g+hgh^2+gh^2g+hg^2h=0.
\end{split}
\end{equation}
The element  $g^3h$ must coincide with  one of the following elements:

Case 1. Let $g^3h=(hg)^2$. Clearly, $h=g(hg)^2$ and  $h^2=(gh)^3$, so $|gh|=2$, a contradiction.

Case 2. Let $g^3h=ghg^2$. Then  $(h, g^2)=1$ and (\ref{E:3}) can be rewrite as
\[
\begin{split}
(gh)^2&+(hg)^2+h^2gh+gh^3+h^3g+hgh^2+gh^2g+g^2h^2=0.
\end{split}
\]
Then $(g,h^2)=1$ and  $(gh)^2=(hg)^2$.
It follows that
\[
H=\gp{\quad g,h \mid g^4=h^4=(g, h^2)=(g^2, h)=1, (gh)^2=(hg)^2\quad}\cong G_{32}^2.
\]
Case 3. Let  $g^3h=h^2g^2$. Then $h=gh^2g^2$ and $hg=gh^2g^3$, so $(hg)^2=1$ which is impossible.

Case 4. Let $g^3h=h^2gh$ or $g^3h=gh^3$.  Then $g^2=h^2$  and  by (\ref{E:3}) $(gh)^2=(hg)^2$. It follows that
\[
H=\gp{ \quad g,h \mid g^4=h^4=1, g^2=h^2,(gh)^2=(hg)^2 \quad }\cong G_{16}^{4}.
\]

Case 5. Let  $g^3h=hgh^2$. Then $|gh|=2$, a contradiction.

Case 6. Let  $g^3h=gh^2g$. Then  $h=g^2h^2g$, so $gh=g^{-1}(h^2)g$ and $2=|h^2|=|gh|$, a contradiction.

Case 7.  Let $g^3h=h^3g$. Then  $gh^3=hg^3$ and by  (\ref{E:3}) we get that
\begin{equation}\label{E:4}
\begin{split}
ghgh+hghg&+\underline{ghg^2}+g^2h^2+h^2g^2+\\
&+g^2hg+h^2gh+gh^2g+hgh^2+hg^2h=0.
\end{split}
\end{equation}
It is easy to check that  $ghg^2\in \{ (hg)^2, h^2gh,  hgh^2, hg^2h\}$. We consider each case separately.

Case 7.1. Let $ghg^2=hghg$. Then  from   (\ref{E:4}) follows that
\begin{equation}\label{E:5}
g^2h^2+h^2g^2+hgh^2+\underline{h^2gh}+gh^2g+hg^2h=0.
\end{equation}
It is easy to check that only possible cases are $h^2gh\in\{ g^2h^2, gh^2g \}$.

If  $h^2gh=g^2h^2$  then  $g^2h=h^2g$  and
$h^3g=g^3h=g(g^2h)=g(h^2g)$, so  $h=g$, a contradiction.

If  $h^2gh=gh^2g$  then  $h(hgh)=(ghg)g^3hg$, then $h(ghg)=(hgh)g^3hg$, so $1=ghg^2$, a contradiction.

Case 7.2. Let $ghg^2=h^2gh$.  Multiplying it on the left side by $g^2$ and on the right size by $gh$ we obtain that
$1=(g^3h)^2=g^2h^2ghgh$. Since $|gh|=4$, this yields that $g^2h^2=ghgh$ and $(g,h)=1$, a contradiction.

Case 7.3. Let $ghg^2=hg^2h$.  Multiplying it on the left side by $g^2$ and on the right size by $gh$ we obtain that
$1=(g^3h)^2=g^2hg^2hgh$. Since $|gh|=4$, this yields that $g^2hg=ghgh$ and $ghg=hgh$. Clearly  $hghg=ghg^2=hg^2h$, so $(g,h)=1$, a contradiction.

Case 7.4. Let $ghg^2=hgh^2$.   Then
\[
\begin{split}
H=\gp{ \quad g,h\mid g^4&=h^4=1,  (g^3h)^2=1, ghg^2=hgh^2 \quad}\\
&\cong \big((C_4 \times  C_2) \rtimes  C_2\big) \rtimes  C_2\cong G_{32}^{6}.
\end{split}
\]
Put $w=1+g(1+h)\in V(FH)$. It is easy to check that
\[
w^2=1+g^2+(gh)^2+g^2h+ghg.
\]
The powers of $w$ and their orders are easy to calculate using the  package  LAGUNA  of the computational algebra system GAP \cite{GAP}.

However we  assume  $w^4=1$.  By a straightforward calculation
\[
\begin{split}
w^4=h+g^2+(gh)^4&+(g^3h)g+g(ghg^2)+(ghg^2)hg\\
&+g^2(hg)^2+g(hg^3)+(g^2h)^2\\
=h+g^2+(gh)^4&+h^3g^2+g(hgh^2)+(hgh^2)hg\\
&+g^2(hg)^2+g^2h^3+(g^2h)^2=1.\\
\end{split}
\]
Comparing the element $h$ with other elements we have  $h=g^2(hg)^2$. This yields that
$hg=g^2h(ghg^2)=g^2h(hgh^2)$ and $g^{-1}hg=(gh^2)^2$. However $4=|h|>|(gh^2)^2|\leq 2$ because $\exp(G)=4$, a contradiction.
Consequently $exp\big(V(FH)\big)>4$, which is impossible.
\endproof

\begin{corollary}\label{C:1}
If  $\exp\big(V(FG)\big)=4$,   then $G'\leq \Phi(G)\leq  \zeta(G)$, $\Phi(G)$ is elementary abelian and $G$ has nilpotency class $2$.
\end{corollary}
\proof
Let $H=\gp{a,b\in G \mid c=(a,b)\not=1}$.    Clearly $c=g_1^2g_2^2\cdots g_n^2$ for some $g_1,\ldots,g_n\in G$ (see  Theorem 10.4.3  in \cite{Hall_book}, p.178).

Using induction on $n\geq 1$, let us prove that $(c,x)=(g_1^2\cdots g_n^2,x)=1$ for any $x\in G$.

Base of induction: $n=1$.  Then
\quad
$(g_1^2,x)=(g_1,x)(g_1,x,g_1)(g_1,x)=(g_1,x)^2=1$  by  (\ref{E:1}) and Lemma \ref{L:1}.

Put $w=g_1^2g_2^2\cdots g_{n-1}^2$. Using the same arguments, (\ref{E:1}) and Lemma \ref{L:1}
\[
(wg_n^2,x)=(w,x)(w,x,g_n^2)(g_n^2,x)=(g_n^2,x)=(g_n,x)^2=1.
\]\endproof

\begin{lemma}\label{L:2}
Let $G$ be a finite $2$-group, such that its  Frattini subgroup  $\Phi(G)$ is   central elementary abelian, $G'\leq \Phi(G)$ and  $|G'|\leq 4$. If  $char(F)=2$  then  the exponent of  $V(FG)$ is equal $4$.
\end{lemma}

\proof
Let $G=g_1\Phi(G)\cup\cdots\cup g_m \Phi(G)$, where $g_1=1$.  Then any  $u\in V(FG)$ can be written as
$u=\sum_{i=1}^m g_iu_i$, where  $u_1,\ldots, u_m\in F\Phi(G)$. By  Brauer's  lemma (\cite{Bovdi_book}, Proposition 3.1, p.17)
\quad $u^2=\sum_{i=1}^m g_i^2u_i^2 +\sum_{1<i<j}^m [g_i,g_j]u_iu_j$. where $[g_i,g_j]=g_ig_j-g_jg_i\in FG$. The element $\sum_{i=1}^m g_i^2u_i^2$ is a central involution, so
\[
u^4 =1+\Big(\sum_{1<i<j}^m [g_i,g_j]u_iu_j\Big)^2.
\]
Since
$[g_i,g_j]=g_ig_j\big(1-(g_j,g_i)\big)$ and $1-(g_j,g_i)$ is a central  nilpotent element of index 2, by Brauer's  lemma (\cite{Bovdi_book}, Proposition 3.1, p.17)  we have that
\[
\begin{split}
z=\Big(\sum_{1<i<j}^m g_ig_j\cdot \big(1-&(g_j,g_i)\big)u_iu_j\Big)^2=\sum_{1<i<j,1<k<l}^m g_ig_jg_kg_l\times \\ &\times \big(1-(g_kg_l,g_ig_j)\big)\big(1-(g_j,g_i)\big)\big(1-(g_l,g_k)\big)u_iu_ju_ku_l.
\end{split}
\]
Suppose that $G'\cong C_2\times C_2$. It is easy to check   that $1-(g_j,g_i)\in \omega(FG')$ and $\omega(FG')^3=0$. Consequently, $z=0$ and
$\exp\big(V(FG)\big)=4$. \endproof

\begin{lemma}\label{L:3}
Let  $G=H\times\gp{a}$ and  $|a|$ divides $4$. If  $exp (V(FH))=4$, then
$exp(V(FG))=4$, too.
\end{lemma}

\proof
First assume that $|a|=2$. Any $x\in V(FG)$ has the form
\[
x=\sum_{g\in H}\alpha_gg+a\sum_{h\in H}\beta_hh
\]
and $g^2, h^2\in \Phi(G)$ are central by Corollary \ref{C:1}. Clearly, the order of the unit
$y=\sum_{g\in H}\alpha_gg+\sum_{h\in H}\beta_hh$
divides $4$ and
\[
y^2=\sum_{g^2\in H}\alpha_g^2g^2+\sum_{h^2\in H}\beta_h^2h^2+\sum_{g,h\in H}\alpha_g\beta_h[g,h].
\]
The unit\quad  $\sum_{g^2\in H}\alpha_g^2g^2+\sum_{h^2\in H}\beta_h^2h^2$ is central and its   order divide $2$, so
\[
y^4=1+\Big([\sum_{g\in H}\alpha_gg,\sum_{h\in H}\beta_hh]\Big)^2=1
\]
and $\big([\sum_{g\in H}\alpha_gg,\sum_{h\in H}\beta_hh]\big)^2=0.$
It follows that  $|x|$ also divides $4$.

Finally  let $|a|=4$ and $L=H\times\gp{a^2}$. Then $exp\big(V(FL)\big)=4$
and any $x\in V(FG)$ has the form $x=\sum_{g\in L}\alpha_gg+a\sum_{h\in L}\beta_hh$.
Repeating  the previous  argument it is easy  to see  that $|x|$ divides $4$.
\endproof

\begin{lemma}\label{L:4}
Let $char(F)=2$ and
let $G$ be a finite $2$-group, such that  $G'$ is  central elementary abelian. If  $|G'|\geq 8$,  then   $\exp\big(V(FG)\big)>4$.
\end{lemma}

\proof
If\quad  $x=\sum_{g\in G}\alpha_gg\in V(FG)$, then\quad  $x^2=\sum_{g\in G}\alpha_g^2g^2+w_2$, where
$w_2\in \gp{g_ig_j(1+(g_i,g_j)) \mid i,j\in\mathbb{N}}_F$. Since $g^2\in \Phi(G)\subseteq \zeta(G)$, we have
\[
x^4=\sum_{g\in G}\alpha_g^4g^4+w_2^2= \sum_{g\in G}\alpha_g^4+w_2^2.
\]
Using the equalities  (\ref{E:1}) and the fact that $|G'|\geq 8$, we have that
\[
\begin{split}
w_2^2\in \gp{g_ig_jg_kg_l\big(1+(g_ig_j,g_kg_l)\big)\big(1+(g_i,g_j)\big)&\big(1+(g_l,g_k)\big) \mid i,j,k,l}_F
\not=0\\
\end{split}
\]
 so $x^4\not=1$. \endproof
\proof[Proof of the Theorem] Follows from   Lemmas \ref{L:1}--\ref{L:4}.\endproof


{\small
{\em Authors' addresses}:

{\em V. Bovdi, M. Salim}, Department of Math. Sciences, UAE University, Al-Ain,\newline United Arab Emirates

 e-mail: \texttt{\{v.bodi, msalim\}@uaeu.ac.ae}

}

\end{document}